\documentclass[]{interact}
\usepackage{epstopdf}
\usepackage[caption=false]{subfig}

\usepackage[mathscr]{eucal}
\usepackage[cp1251]{inputenc}
\usepackage[T2A]{fontenc}
\usepackage{amsfonts,amssymb,amsmath,latexsym}

\usepackage[numbers,sort&compress]{natbib}
\bibpunct[, ]{[}{]}{,}{n}{,}{,}
\renewcommand\bibfont{\fontsize{9}{13}\selectfont}
\makeatletter
\def\NAT@def@citea{\def\@citea{\NAT@separator}}
\makeatother

\theoremstyle{plain}
\newtheorem{theorem}{Theorem}[section]

\theoremstyle{definition}

\theoremstyle{remark}

\begin{document}

\title{Singular Behavior of Harmonic Maps Near Corners}

\author{
\name{S.~I. Bezrodnykh\thanks{CONTACT S.~I. Bezrodnykh. Email:
sbezrodnykh@yandex.ru, sbezrodnykh@mail.ru}, V.~I.
Vlasov\thanks{CONTACT V.~I. Vlasov. Email: vlasov@ccas.ru}}
\affil{\textsuperscript{a} Dorodnicyn Computing Centre, FRC CSC
RAS, 40 Vavilova st., Moscow, Russia, 119333\\
 Sternberg
Astronomical Institute, Lomonosov Moscow State University, 13
Universitetskii prosp., Moscow, Russia, 119992\\
Peoples Frienship University of Russia (RUDN University), 6
Miklukho-Maklaya st., Moscow, Russia, 117198.}
\affil{\textsuperscript{b} Dorodnicyn Computing Centre, FRC CSC
RAS, 40 Vavilova st., Moscow, Russia, 119333\\
 Peoples Frienship
University of Russia (RUDN University), 6 Miklukho-Maklaya
 st., Moscow, Russia, 117198.} }

\maketitle

\begin{abstract}
    For a harmonic map $\mathcal{F}:\mathcal{Z} {\buildrel {\,harm\,}
\over\longrightarrow} \mathcal{W}$ transforming the contour of a
corner of the boundary $\partial\mathcal{Z}$ into a rectilinear
segment of the boundary $\partial\mathcal{W}$, the behavior near
the vertex of the specified corner is investigated.
 The behavior of the inverse map $\mathcal{F}^{-1}:\mathcal{W}
\longrightarrow \mathcal{Z}$ near the
preimage of the vertex is investigated as well.
 In particular, we prove that if $\varphi$ is the value of the exit angle from
 the vertex of the reentrant corner for a smooth curve $\mathcal{L}$ and $\theta$
 is the value of the exit angle from the vertex image for the image
$\mathcal{F} (\mathcal{L})$ of the specified curve, then the
dependence of $\theta$ on $\varphi$ is described by a
discontinuous function.
 Thus, such a behavior of the harmonic map qualitatively differs from the behavior
 of the corresponding conformal map: for the latter one, the dependence $\theta
(\varphi)$ is described by a linear function.
\end{abstract}


\begin{amscode}
AMS 2010 mathematical subject classification: Primary 35J40; Secondary: 35J58, 35J30, 30G20
\end{amscode}

\section{Introduction}

We say that a homeomorphic map of domains $\mathcal{Z}$ and
$\mathcal{W}$ located in the complex planes  $z = x + i y$ and $w
= u + i v$ respectively is \emph{harmonic} (see, e.\,g.,
\cite{Hamilton1975}--\cite{Bshouty2005}) and denote it by the
symbol
\begin{equation}
\mathcal{F}: \mathcal{Z} {\,\buildrel {harm\,\,\,}
\over\longrightarrow\,} \mathcal{W} \label{BV-1}
\end{equation}
 if the corresponding complex-valued
 function $\mathcal{F} (z) = u (x, y) + i v (x, y)$ is harmonic.
 Note that, unlike conformal maps, we do not require the scalar harmonic
 functions $u$ and $v$ to obey the Cauchy--Riemann
 conditions.

If $\mathcal{Z}$ and $\mathcal{W}$ are Jordan domains, then  map
(\ref{BV-1}) is frequently constructed  (see, e.\,g.,
\cite{Hamilton1975}--\cite{Choquet1945}) as a harmonic extension
of the given homeomorphism  $\mathcal{B}:
\partial\mathcal{Z} {\,\buildrel {Hom\,\,\,}
\over\longrightarrow\,}
\partial\mathcal{W}$ of their boundaries into the domain $\mathcal{Z}$,
 i.\,e., as a solution of the following Dirichlet problem for the
 Laplace equation:
\begin{equation}
\Delta\, \mathcal{F} (z)\, =\, 0, \quad z \in \mathcal{Z};\qquad
\mathcal{F} (z)\, =\, \mathcal{B} (z), \quad z \in  \partial  \mathcal{Z}.
\label{BV-2}
\end{equation}
    However, it is not guaranteed that the map constructed that way is one-sheeted;
 this is easily confirmed by the corresponding counterexamples
(see, e.\,g., \cite{BezrodnykhVlasov-2012}).
 Note that, for conformal maps (unlike harmonic ones), the homeomorphism of the boundaries
 of Jordan domains implies the homeomorphism of the (closed)
 domains themselves (see the Caratheodory theorem in \cite{Caratheodory1913}). 
 For the map $\mathcal{F} (z)$ defined by (\ref{BV-2}), a sufficient condition to be
 one-sheeted is the convexity of the domain $\mathcal{Z}$ (see the
  Rad\'o--Kneser--Choquet
theorem in \cite{Duren2004},\cite{Rado1926}--\cite{Choquet1945}).

Harmonic maps are broadly investigated (see
\cite{Hamilton1975}--\cite{Bshouty2005},\cite{BezrodnykhVlasov-2012},\cite{Renelt1988}--\cite{Alessandrini2009}).
 The interest to this direction is caused by its great theoretical importance.
 In particular, harmonic maps play a substantial role in differential geometry and
 the theory of minimal surfaces (see \cite{Bers1951},\cite{Nitsche1958}),
 in the geometric function theory  (see \cite{Clunie1984}--\cite{Bagapsh2017}),
 and in the geometric theory of quasiconformal maps
 (see \cite{Grotzsch1928}--\cite{Vuorinen1997}),
 including the extension problem of the Riemann theorem for those
 maps (see \cite{Duren2004},\cite{Bers1953}--\cite{DurenKhavinson1997}).

On the other hand, harmonic maps have numerous applications.
  In particular, they are applied for the constructing of     computational
  meshes in complicated domains  $\mathcal{Z}$, based on the
  Winslow approach (see \cite{Winslow1967}-\cite{Azarenok2014}).
  Under such an approach, the desired mesh  (in a Jordan domain) is constructed as follows:
  we take the Cartesian mesh natural for the square $\Pi := \{u \in (- 1/2, 1/2), v \in (0,
1)\}$ and transform it to the domain $\mathcal{Z}$ by means of the
map
\[
\mathcal{F}^{-1} (w)\, =\, x (u,\, v) + i y (u,\, v) :\,\,
\Pi {{\,\buildrel {\,\,harm^{-1}} \over\longrightarrow\,}} \mathcal{Z}\,,
\]
 which is inverse to a harmonic map.
 To find the specified map, we use numerical methods to solve the following Dirichlet problem
 in the square  $\Pi$ for the ``inverse vector'' Laplace equation:
\begin{equation}
\Delta^{-1}\, \mathcal{F}^{-1} (w)\, =\, 0, \quad w \in \Pi;\qquad
\mathcal{F}^{-1} (w)\, =\, \mathcal{B}^{-1} (w), \quad w \in  \partial  \Pi.
\label{BV-3}
\end{equation}
 The expanded form of the said equation is the following quasilinear system
 of equations with respect to the component pair $x (u,\, v),\, y (u,\, v)$:
\[
\left\{\begin{array}{c}
\vspace*{2mm}
(x_v^2  +  y_v^2) x_{uu}  -2 (x_u x_v + y_u y_v) x_{uv} +
(x_u^2  +  y_u^2) x_{vv}  = 0,\\
(x_v^2  +  y_v^2) y_{uu} -2 (x_u x_v\ + y_u y_v) y_{uv} + (x_u^2
+  y_u^2) y_{vv} = 0;
\end{array}
\right.
\]
 here, subscripts denote the differentiating with respect
  to the corresponding variables (that expanded form can be found, e.\,g.,
  in \cite{Winslow1967},\cite{KnSt1993},\cite{Ivanenko1997}--\cite{Liseikin1999}).

 To pose the boundary-value
 condition for problem (\ref{BV-3}), a homeomorphism  $\mathcal{B}^{-1}$
 between the boundaries $\partial\Pi$ and
$\partial\mathcal{Z}$ is taken such that the following property
takes place: if a point $w$ moves along
 the square side with a constant velocity, then its image
 $z = \mathcal{F}^{-1} (w)$ moves along the boundary
$\partial\mathcal{Z}$ with a constant velocity as well.
 This requirement is imposed to ensure the desired quality of the computational mesh
 (see \cite{KnSt1993},\cite{Ivanenko1997},\cite{ThSW1999}).

   It is convenient to use the map $\mathcal{F}^{-1} (w)$ to construct the mesh in
    a complicated domain $\mathcal{Z}$ because a     computational scheme
    for the resolving of problem (\ref{BV-4}) is rather easily constructed on the Cartesian
    mesh of the square $\Pi$ (see,
 e.\,g., \cite{Godunov1972},\cite{Roache1985},\cite{KnSt1993}--\cite{Liseikin1999});
 on the other hand, the  homeomorphism of the map is guaranteed by the
   Rad\'o--Kneser--Choquet theorem
 (see \cite{Duren2004},\cite{Rado1926}--\cite{Choquet1945})
 because the square is a convex domain.

 The described approach is broadly propagated (see \cite{Winslow1967}-\cite{Azarenok2014}),
 but it occurs that its practical applications face difficulties
 for the numerically obtained map: it even might lose the
 homeomorphism (see, e.\,g., \cite{Roache1985},\cite{KnSt1993}-\cite{Ivanenko1997}),
 contradicting the above theoretical aspects.
 In particular, such a difficulty arises if the domain
$\mathcal{Z}$ has reentrant corners, i.\,e., with angles exceeding $\pi$
(see \cite{Ivanenko1997},\cite{Azarenok2014}).
 To overcome that obstacle, one has (primarily) to understand its nature.
 To do that, one has to consider harmonic maps $\mathcal{F}$ transforming
 the contour of a corner $\partial\mathcal{Z}$ of the boundary
 into a rectilinear segment of the boundary $\partial\Pi$
 and investigate its behavior near the corner vertex
(this is provided in Sec.  4) as well as the behavior of the
inverse map near the  vertex preimage (see Sec.  3).
 To obtain those results, we use the asymptotic behavior of the
 map $\mathcal{F}$  near the corner vertex, found in Sec.  2. Also, we compare the harmonic
  map $\mathcal{F}$ with the corresponding conformal map (see the final part of Sec.  4).

\begin{figure}
\centering
\resizebox*{15cm}{!}{\includegraphics{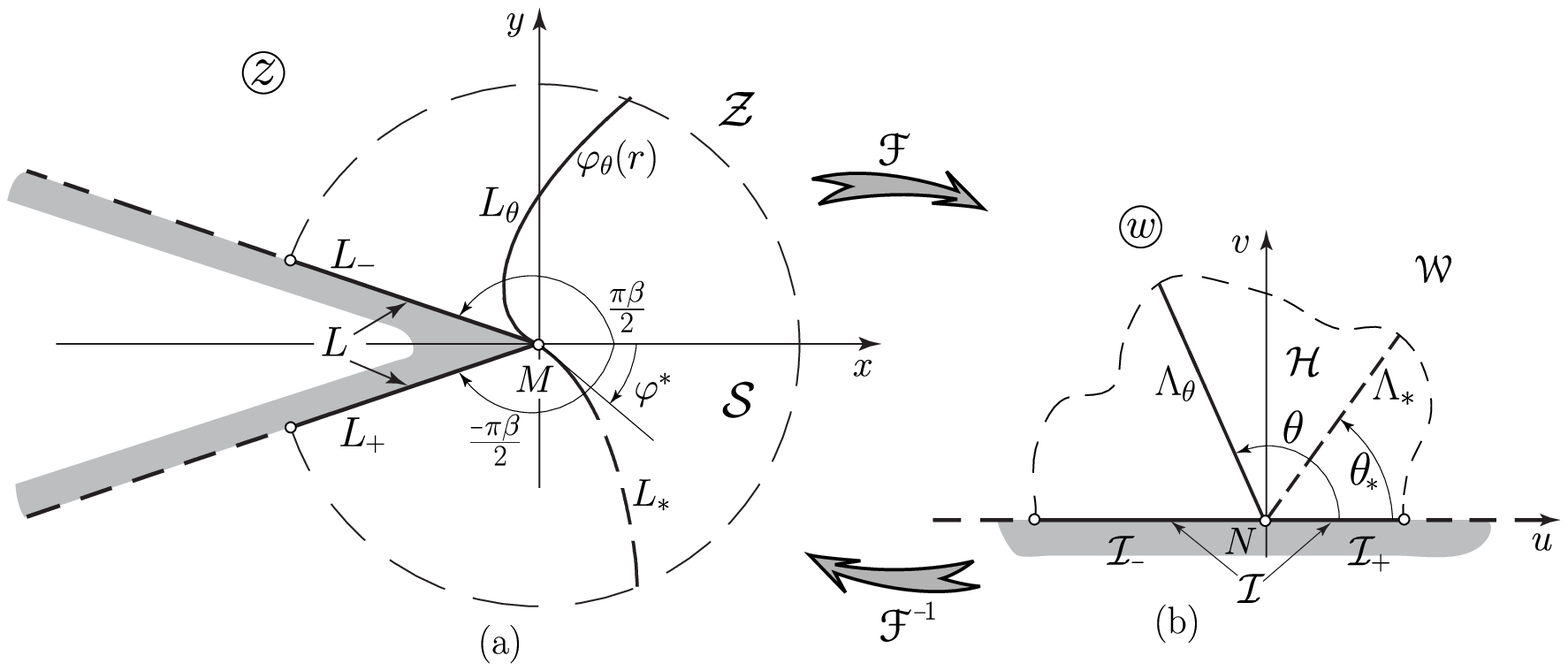}}
\caption{} \label{fig1}
\end{figure}

\section{Harmonic maps near corner vertices: asymptotic behavior}

Let $\mathcal{Z}$ be a domain of the complex plane $z = r\,e^{ i
\varphi}$. Let its boundary $\partial\mathcal{Z}$ contain an angle
of value $\pi \beta$ (measured with respect to the domain), $\beta
\in (0, 2)$, such that its vertex is the co-ordinate origin and its sides are
$L_{\pm} =\bigl\{r \in (0, R), \varphi = \mp \pi\beta / 2\bigr\}$.
 Then the domain $\mathcal{Z}$ contains a sector
\begin{equation}
\mathcal{S} = \Big\{r \in (0, R),
\varphi \in \Big(- \frac{\pi\beta}{2},\, \frac{\pi\beta}{2}\Big)\Big\}
\label{BV-4}
\end{equation}
 of radius $R$, adjoining the boundary $\partial\mathcal{Z}$
 along the corner contour $L = L_- \cup \{0\} \cup L_+$
  (see Fig. 1a).
  Further, let $\mathcal{W}$ be a domain of the complex plane $w = \rho e^{i\theta}$
  such that its boundary contains a segment $\mathcal{I} \ni \{0\}$
  of the real axis
$\mathbb{R}$ (see Fig. 1b) and the domain  $\mathcal{W}$ itself
 contains a subdomain  $\mathcal{H}$ of the upper half-plane
$\mathbb{H}$, adjoining the segment $\mathcal{I}$.
 Finally, let a harmonic function $w = \mathcal{F} (z)$ maps
$\mathcal{Z}$ onto $\mathcal{W}$ such that the mapping of the
sector $\mathcal{S}$ into $\mathcal{H}$ is one-sheeted
 and the homeomorphism $w = \mathcal{B} (z)$ between $L$ and
$\mathcal{I}$ (see Fig. 1) acts as follows:
\begin{equation}
u (z) = - \sigma_{_-}\, r, \;\, z \in L_-,\qquad u (z) =
\sigma_{_+}\, r, \;\, z \in L_+, \qquad v (z) = 0,\;\, z \in L,
\label{BV-5}
\end{equation}
 where $\sigma_{_-}$ and $\sigma_{_+}$ are positive numbers.
 Note that the required rule corresponds to the coordination requirement
 for the ``movement velocities'' of the point $z$ on $L$ and its image $w$
 on $\mathcal{I}$, imposed in  Sec.  1.
 Homeomorphism (\ref{BV-5}) preserves the orientation of the boundaries:
 if the contour $L$ of the corner is passed through its points $z$ to leave
 the domain $\mathcal{Z}$ on the left, then the contour  $\mathcal{I}$ is passed
  through points $w = \mathcal{B} (z)$ such that $\mathcal{H}$ is left on the left
  as well. This implies that $\mathcal{B} (L_+)=:
\mathcal{I}_+ \subset \mathbb{R}_+$, $\mathcal{B} (L_-)=:
\mathcal{I}_- \subset \mathbb{R}_-$, and $\mathcal{B} (0) = 0$
(see Fig.~1).

Introduce the linear function
\begin{equation}
Q (z):= - \mu\, r \sin (\varphi -\varphi^*)
\label{BV-6}
\end{equation}
 such that its parameters $\mu$ and $\varphi^*$ are defined by the angle value $\pi \beta$
 and the ``velocities'' $\sigma_-$ and $\sigma_+$ of homeomorphism (\ref{BV-5}) as follows:
\begin{equation}
\mu := \,|\sin \pi\beta|^{-1}\, {\sqrt{\sigma_{_+}^2\, +\,
\sigma_{_-}^2 \,+ \,2 \sigma_{_+} \sigma_{_-} \cos \pi
\beta}\,}\;~ \textrm{and} ~\;\varphi^* := \arctan
\biggl(\frac{\sigma_{_+} - \sigma_{_-}}{\sigma_{_+} + \sigma_{_-}}
\tan \frac{\pi \beta}{2} \biggr). \label{BV-6a}
\end{equation}
 The imposed positivity of $\sigma_+$ and $\sigma_-$ implies that
$\mu > 0$ and $\varphi^* \in \bigl(-\pi\beta/2, \pi\beta/2\bigr)$.

 One can verify that function (\ref{BV-6}) satisfies Condition (\ref{BV-5})
 on the contour $\Gamma$ of the corner, i.\,e.,
  \[ Q (r\,e^{ i
\varphi}) \,=\, \left\{
\begin{array}{rl}
  \sigma_{_+}\, r,   &\varphi = - \pi\beta/2, \\
- \sigma_{_-}\, r,   &\varphi = \pi\beta/2.
\end{array}\right.
\]
 Therefore, the function $\mathcal{F} - Q$, which is harmonic in sector (\ref{BV-4}),
 vanishes on $L$ and, therefore, on the set $\mathcal{S} \cup L$,
 it is represented by a series
\begin{equation}
\mathcal{F} (z) - Q (z) = \sum_{n = 1}^\infty\,(a_n + i b_n)\,\mathrm{Im} \Bigl(
i\, z^{1 / \beta} \Bigr)^n\,,
\label{BV-7}
\end{equation}
 where $a_k$ and $b_k$ are real numbers and series \eqref{BV-7} converges on
 the set $\mathcal{S} \cup L$. Since $\mathcal{F}$ is a one-sheeted
  map preserving its orientation, it follows that
\begin{equation}
a_1 \not= 0 \qquad\textrm{and} \qquad b_1 > 0. 
\label{BV-8}
\end{equation}
  Move the function $Q$ defined by relation (\ref{BV-6}) to the right-hand
  part of (\ref{BV-7}) and truncate its series after the third term.
  This yields the asymptotic behavior of the considered harmonic map
  near the corner vertex  $z = 0$, uniform with respect to $\varphi
$ from $ [-\pi \beta / 2,\,\pi \beta / 2]$.  For $\beta \in (0,
1)$, the specified asymptotic behavior is as follows:
\begin{equation}
\begin{split}
\mathcal{F} (r e^{ i \varphi})  = - \mu  r \sin (\varphi
-\varphi^*) + (a_1 + i b_1)  r^{1 /\beta} \cos
\dfrac{\varphi}{\beta}\\
 + (a_2 + i b_2)  r^{2 /\beta} \sin
\dfrac{2\varphi}{\beta}
+ \mathcal{O} \bigl(r^{3 /\beta} \bigr), \quad r \to 0.
\end{split}
 \label{BV-9}
\end{equation}
 If $\beta \in (1, 2),$ then the order of the second term is greater than the order of the
 first one; therefore, they are to interchange their places.

\section{The map inverse to the harmonic one:\\
 behavior near the preimage of the corner vertex}

 To investigate the behavior of the map
\[
z = \mathcal{F}^{-1} (w): \mathcal{H} {\,\buildrel {harm^{-1}}
\over\longrightarrow\,} \mathcal{S}
\]
 inverse to the map $\mathcal{F}$ near the point $w =
0$, which is the image $\mathcal{F} (0)$ of the corner vertex (see
Sec.  2), we introduce the curves $L_\theta$ on the plane $z$ as
follows: they are the images  (under the inverse map) of the rays
\[
\Lambda_\theta := \{\mathcal{H} \ni w: |w| > 0, \arg w = \theta\}
\]
 starting from the point $w = 0$ under the angle $\theta$, i.\,e., $L_\theta
:= \mathcal{F}^{-1} (\Lambda_\theta)$ (see Fig. 1); recall that
$(\rho, \theta)$ are the polar coordinates on the plane  $w$.

Let $\varphi = \varphi_\theta (r)$ be the equation of the curve
 $L_\theta$ in the polar coordinates $(r, \varphi)$. Let us find the asymptotic behavior
 of the function $\varphi_\theta (r)$ as  $r \to 0$. Since the curve $L_\theta$
 satisfies the relation ${\rm Im}
\bigl[\mathcal{F} (L_\theta) \, e^{- i \theta}\bigr] = 0$ by
definition, it follows that the function $\varphi_\theta (r)$
 satisfies the relation
\[
{\rm Im} \bigl[\mathcal{F} \bigl(r e^{ i \varphi_\theta (r)}\bigr)\,
e^{- i \theta}\bigr] = 0.
\]
  In that last relation, take into account the asymptotic behavior of
  the map $\mathcal{F}$, given by (\ref{BV-9}).
  We obtain the following (asymptotic) equation for the desired function $\varphi_\theta (r)$:
\begin{equation}
\begin{split}
\mu\, r \sin \bigl(\varphi_\theta (r) -\varphi^*\bigl) \, =\,
(a_1 -  b_1\,\cot \theta)\, r^{1/\beta} \cos \frac{\varphi_\theta (r)}{\beta}\\
+\, (a_2 -  b_2\,\cot \theta)\, r^{2/\beta} \sin \frac{2\, \varphi_\theta (r)}{\beta}\,+\,
\mathcal{O} \bigl(r^{3 /\beta} \bigr),\qquad r \to 0.
\label{BV-10}
\end{split}
\end{equation}
    First, consider the special case where the expression in
    brackets at the right-hand
    part of (\ref{BV-10}) vanishes, i.\,e.,
 the following relation is satisfied:
\begin{equation}
a_1 -  b_1\,\cot \theta = 0.
\label{BV-10a}
\end{equation}
   Denote the angle corresponding to that relation as follows:
\begin{equation}
\theta^* = \arctan\, \frac{b_1}{a_1}. \label{BV-10b}
\end{equation}
 The curve $L_{\theta}$ corresponding to that case is denoted by $L^*$.
 The equation for that curve is denoted by  $\varphi = \varphi^* (r)$.
 Its geometric properties are different from the ones of the curve $L_\theta$
 for $\theta \not=\theta_*$.

Substitute (\ref{BV-10a}) in relation (\ref{BV-10}) and divide
both parts of the obtained relation by $\mu r$; this yields the
following (asymptotic) equation for $\varphi^* (r)$:
\[ 
\sin \bigl(\varphi^* (r) -\varphi^*\bigl) \, =\,
\frac{(a_2 -  b_2\,\cot \theta)}{\mu}\, r^{2/\beta - 1}
\sin \frac{2\, \varphi^* (r)}{\beta}\,+\,
\mathcal{O} \bigl(r^{3 /\beta - 1} \bigr),\qquad r \to 0.
\]
    Taking into account that its right-hand
 part tends to zero as $r \to 0$ provided that $\beta \in (0, 2)$,
 we obtain the desired asymptotic behavior in the form
\begin{equation}
\varphi^* (r)\,=\, \varphi^*\,+\, E_1^*\, r^{2 /\beta - 1}\, +\, o \bigl(r^{2 /\beta - 1}
\bigr),\quad r \to 0, \qquad \theta = \theta^*,\qquad \beta \in (0, 1)\,,
\label{BV-10c}
\end{equation}
 where
\[
E_1^*\,=\, \mu^{- 1}\, (a_2 -  b_2\,\cot \theta) \sin \frac{2 \varphi^*}{\beta}\,.
\]
  Now, let $\theta \not=\,\theta^*$. 
 Then the following two cases are to be selected:
\begin{equation}
({\rm I})\, \beta \in (0, 1)\qquad \textrm{and} \qquad ({\rm
II})\, \beta \in (1, 2). \label{BV-11}
\end{equation}
 In the second case, we say that the corner is reentrant.
 This case is the most interesting from the viewpoint of applications and the difficulty
 treated in Sec. 1 arises in this case. However, we start our
 investigation from the first case.

({\rm I}) Let $\beta \in (0, 1)$. Then, dividing both parts of
(\ref{BV-10}) by $\mu\,r$ we obtain the equation
\begin{equation}
\begin{split}
\sin \bigl(\varphi_\theta (r) -\varphi^*\bigl) \, =\,
\mu^{-1}(a_1 -  b_1\, \cot \theta)\, r^{1/\beta-1}
\cos \frac{\varphi_\theta (r)}{\beta}\\
+\, \mu^{-1} (a_2 -  b_2\,\cot \theta)\, r^{2/\beta-1} \sin
\frac{2\, \varphi_\theta (r)}{\beta}\,+\, \mathcal{O} \bigl(r^{3
/\beta-1} \bigr),\qquad r \to 0.
\end{split}
\label{BV-11a}
\end{equation}              
    Its right-hand
 part tends to zero as $r \to 0$. Taking this fact into account, we find the
 following
 asymptotic behavior for the considered case (I):
\begin{equation}
\varphi_\theta (r)\,=\, \varphi^*\,+\, E_1\, r^{1 /\beta - 1}\, +\, o \bigl(r^{1 /\beta - 1}
\bigr),\quad r \to 0, \qquad \theta \in (0, \pi),\quad \beta \in (0, 1),
\label{BV-12}
\end{equation}
 where
\[
E_1\,=\, \mu^{- 1}\, (a_1 -  b_1\,\cot \theta) \cos \frac{ \varphi^*}{\beta}\,.
\]
  If $\theta = 0$ or $\theta = \pi$, then we have the following exact relations
  instead of  (\ref{BV-12}), which is an asymptotic one:
\begin{equation}
\varphi_0 (r) = - \frac{\pi \beta}{2}\qquad \textrm{and}\qquad
\varphi_\pi (r) = \frac{\pi \beta}{2}; \label{BV-13}
\end{equation}
 this obviously follows from Conditions (\ref{BV-5}).
 Thus, the dependence $\varphi (\theta)$, where $\theta$ is the exit angle
 of the ray $\Lambda_\theta$ leaving the preimage $w = 0$ of the corner vertex,
 while $\varphi$ is the exit corner
 of the image
$L_\theta = \mathcal{F} (\Lambda_\theta)$ of the ray, leaving the
corner vertex  $z = 0$ itself, is as follows (provided that $\beta
\in (0, 1)$):
\begin{equation}
\beta \in (0, 1):\qquad \varphi (\theta) \,=\, \left\{\begin{array}{ll} -
\pi \beta / 2, &
\theta = 0, \\
 \varphi^*, & \theta \in (0, \pi), \\
\pi \beta / 2,& \theta = \pi.
\end{array}\right.
\label{BV-14}
\end{equation}

\begin{figure}
\centering
\resizebox*{5cm}{!}{\includegraphics{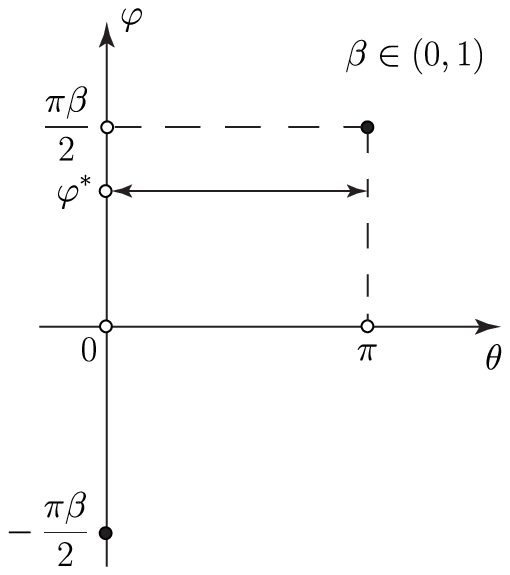}}
\caption{} \label{fig2}
\end{figure}

 Figure 2 displays the graph of dependence (\ref{BV-14}). Figure 3a provides the polar
 mesh
\begin{equation}
\Xi = \Big\{\rho = \frac{n}{5}, n = \overline{1, 5};\; \theta =
\frac{\pi n}{7}, n = \overline{0, 7}\Big\}, \label{BV-14a}
\end{equation}
 while Fig. 3b displays its image $\mathcal{F}^{-1} (\Xi)$ on the plane $z$
 under the inverse map (for the considered case {\rm I}).

\begin{figure}
\centering
\resizebox*{12cm}{!}{\includegraphics{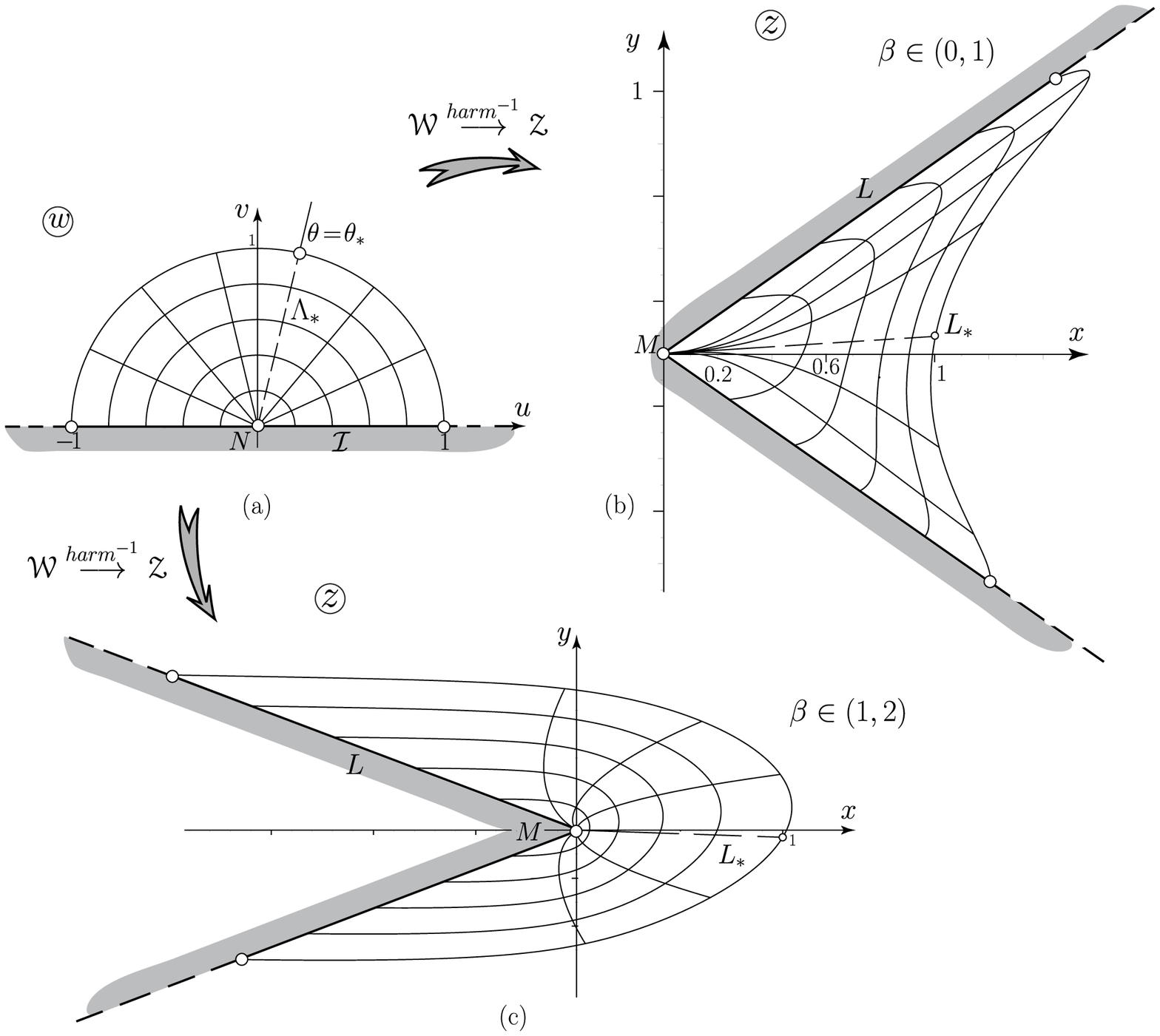}}
\caption{} \label{fig3}
\end{figure}

({\rm II})
 Now, assume that the second condition of (\ref{BV-11}) is satisfied and
  $\theta\not=\theta^*$.
  To obtain the asymptotic behavior of the function $\varphi_\theta (r)$
  describing the curve $L_\theta$, we take into account that
  $(a_1 -  b_1\,\cot \theta)\, r^{1/\beta}$ is different from zero
  and divide Eq. (\ref{BV-10}) by the specified value.
  This allows one to represent the specified equation as follows:
\begin{equation}
\cos \frac{\varphi_{\theta} (r)}{\beta} =
F_1\;r^{1 - 1 /\beta} \sin \bigl(\varphi_{\theta} (r) - \varphi^*\bigr) +
\mathcal{O} \bigl(r^{1 /\beta}\bigr),\quad r \to 0; \qquad
\theta^* \not= \theta \in (0, \pi),
\label{BV-17}
\end{equation} 
 where
\[
F_1  = \mu\, (a_1 -  b_1\,\cot \theta\,)^{-1}\,.
\]
 Since the right-hand
 part of relation (\ref{BV-17}) tends to zero as $r \to 0$, it
 follows that the first term of the asymptotic representation of $\varphi_{\theta} (r)$
 is either $\pi \beta / 2$ or $- \pi \beta / 2$.
 The detailed analysis of Eq. (\ref{BV-17}) and relation (\ref{BV-10c}) yields the following
 result:
\begin{equation}
\varphi_\theta (r) \,=\, \left\{
\begin{array}{cl}
- \frac{\pi \beta}{2} - \Big[\beta\,F_1\,\sin \bigl(\frac{\pi \beta}{2}
+ \varphi^*\bigr) \Big] r^{1 - 1/\beta} + O (r^\gamma)\,, &\quad \theta \in (0, \theta^*),\\
\\
\varphi^*\,+\, E_1^*\, r^{2 /\beta - 1}\, +\,
o \bigl(r^{2 /\beta - 1} \bigr), &\quad \theta = \theta^*\,,
\\ \\
\frac{\pi \beta}{2} - \Big[\beta\,F_1\,
\sin \bigl(\frac{\pi \beta}{2}  - \varphi^*\big)\Big]  r^{1 - 1/\beta}
 + \mathcal{O} (r^\gamma), &\quad \theta \in (\theta^*, \pi),
\end{array}\right.~\; r \to 0,
\label{BV-18}
\end{equation}
 where $\gamma\, =\, \min\big\{1/\beta,\, 2(1 - 1/\beta)\big\}$.
 This implies that if $\beta \in (1, 2)$, then the dependence  $\varphi (\theta)$,
  where $\theta$ is the exit angle
 of the ray $\Lambda_\theta$ leaving the preimage $w = 0$ of the corner vertex,
 while $\varphi$ is the exit angle
 of the image
$L_\theta = \mathcal{F}^{-1} (\Lambda_\theta)$ of the ray, leaving
the corner vertex  $z = 0$ itself, is described by the relation
\begin{equation}
\beta \in (1, 2): \qquad \varphi (\theta) \,=\, \left\{\begin{array}{cl} -
\pi \beta / 2, &\theta \in (0, \theta^*), \\
\quad \varphi^*, &\theta = \theta^*, \\
\pi \beta / 2, &\theta \in (\theta^*, \pi).
\end{array}\right.
\label{BV-18a}
\end{equation}
 Figure 4 displays the graph of dependence (\ref{BV-18a}).
 Figure 3c displays $\mathcal{F}^{-1} (\Xi)$ on the plane $z$,
 i.\,e.,  the image under the inverse map for the polar mesh defined by (\ref{BV-14a})
  (for case {\rm II}).

\begin{figure}
\centering
\resizebox*{5cm}{!}{\includegraphics{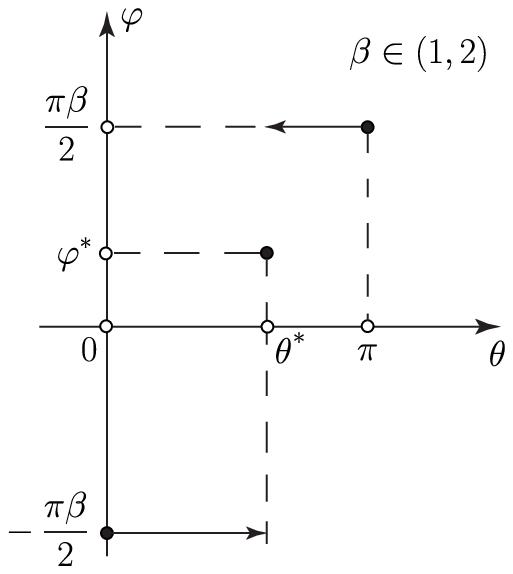}}
\caption{} \label{fig4}
\end{figure}

Combining relation (\ref{BV-18}) with the similar result for case
I, we arrive at the following assertion.

\begin{theorem}\label{theo-1}
The function $\varphi = \varphi_\theta (r)$ describing the curve
$L_\theta := \mathcal{F}^{-1} (\Lambda_\theta)$ in the polar
coordinated obeys the asymptotic relation given by
\eqref{BV-12}-\eqref{BV-13} provided that $\beta \in (0, 1)$ and
 obeys the asymptotic relation given by \eqref{BV-18} provided that $\beta \in (1, 2)$.
\end{theorem}

\section{Behavior of harmonic maps near vertices of corners}

To investigate the behavior of the map $\mathcal{F}: \mathcal{S}
{\,\buildrel {harm\,\,} \over\longrightarrow\,} \mathcal{H}$
 near the vertex $w = 0$ of a corner, introduce curves $l_\varphi$
 on the plane $w =
\rho e^{i\,\theta}$ such that they are the images of the rays
\[
\lambda_\varphi := \{\mathcal{S} \ni z: |z| > 0, \arg z = \varphi\},
\]
starting from the specified  vertex, i.\,e., $l_\varphi :=
\mathcal{F} (\lambda_\varphi)$.
 Similarly to the investigation of the inverse map $\mathcal{F}^{-1}$ (see Sec.  3),
 to investigate the behavior of $\mathcal{F}$ near the point $z = 0$, one has to
 consider cases (\ref{BV-11}) separately.

({\rm I})
 Let $\beta \in (0, 1)$. To find the form of the curve $l_\varphi$
 in the Cartesian coordinates $(u, v)$, one has to find the parametric dependence
 of those coordinates on the coordinates on the plane $z$, i.\,e.,
 to find the functions
$u_\varphi = u (r)$ and $v_\varphi = v (r)$.
 To do that, separate the real and imaginary part in the asymptotic relation given
 by (\ref{BV-9}). We obtain that
\[
u_\varphi (r) =  - \mu r \sin (\varphi - \varphi^*) + a_1
r^{1/\beta} \cos \frac{\varphi}{\beta} + a_2 r^{2/\beta} \sin
\frac{\varphi}{\beta} + O (r^{3/\beta}),\qquad r \to 0,
\]
 and
  \[
v_\varphi (r) = b_1 r^{1/\beta} \cos \frac{\varphi}{\beta} + b_2
r^{2/\beta} \sin \frac{\varphi}{\beta} + O (r^{3/\beta}),\qquad r
\to 0.
\]
 From the former relation, express $r$ via $u$ and $\varphi$. Then substitute it
 in the latter one. We obtain the following asymptotic behavior of the function
  $v_\varphi (u)$ as $\varphi \not=\varphi^*$:
\begin{equation}
\begin{split}
v_\varphi (u) =
\frac{b_1 \cos(\varphi / \beta)}{\big(\mu |\sin (\varphi - \varphi^*\big)|)^{1/\beta}}\,
|u|^{1/\beta} + o (|u|^{1/\beta}),\qquad &|u| \to 0,\\
\beta \in (0, 1), \qquad &\varphi^* \not=\varphi \in
\Big(- \frac{\pi \beta}{2}, \frac{\pi \beta}{2}\Big).
\end{split}
\label{BV-18b}
\end{equation}
 It is substantially simplified if $\varphi =\varphi^*$:
\begin{equation}
v_{\varphi^*} (u) = \frac{b_1}{a_1}\, u + o (u),\quad u\to 0,\qquad
\beta \in (0, 1),\quad \varphi =\varphi^*.
\label{BV-18c}
\end{equation}

({\rm II})
 Let $\beta \in (1, 2)$. Introduce the Cartesian coordinate system
 rotated to the angle $\theta^*$ with respect to the original one.
 Then the new coordinates $(\mathcal{U},\,\mathcal{V})$ are related
 to the original coordinates $(u, v)$ by the relation
\[
\mathcal{U} + i \mathcal{V} = (u + i v)\,e^{ - i \theta^*},
\]
 where $\theta^*$ is the value of the angular coordinate $\theta$ such that the first
 relation of (\ref{BV-10}) is satisfied for  $\theta=\theta^*$.
Let $\mathcal{V} = \mathcal{V}_{\varphi} (\mathcal{U})$ denote the
equation of the curve  $l_\varphi$ in the rotated Cartesian
coordinates $(\mathcal{U}, \mathcal{V})$.

Using asymptotics (\ref{BV-9}) for the map $\mathcal{F}$, we
obtain the desired asymptotic behavior for the function
$\mathcal{V}_{\varphi} (\mathcal{U})$.
 For $\varphi = \varphi^*$, we have the asymptotic relation
\begin{equation}
\mathcal{V}^* (\mathcal{U}) = \frac{2 (a_1 b_2 - a_2 b_1)\,
\tan \frac{\varphi^*}{\beta}}{(a_1^2 + b_1^2)^{3/2}}\,
\mathcal{U}^2 + \mathcal{O} (\mathcal{U}^3), \qquad \mathcal{U} \to 0,
\label{BV-20}
\end{equation}
 where $\mathcal{V}^*$ denotes $\mathcal{V}_{\varphi^*}$. For other values of $\varphi$,
   we have the following asymptotic relation:
\begin{equation}
\mathcal{V}_{\varphi} (\mathcal{U}) = \frac{\mu b_1 \sin (\varphi - \varphi^*)}
{(a_1^2 + b_1^2)^{\frac{\beta + 1}{2}}\, \cos^\beta \frac{\varphi}{\beta}}\,
\mathcal{U}^{\beta} +
 \mathcal{O} (\mathcal{U}^{2 \beta - 1}), \quad \mathcal{U} \to 0, \qquad
 \varphi^* \not= \varphi \in \Big(- \frac{\pi\beta}{2},\, \frac{\pi\beta}{2}\Big).
\label{BV-21}
\end{equation}
    Let $\theta = \theta_\varphi (\rho)$ be the equation of the curve $l_\varphi$
    in the polar coordinates of the plane $w$.
 Represent relations (\ref{BV-11a})-(\ref{BV-12})
  as asymptotic relations
  for $l_\varphi$ in the polar coordinates near the image $w = 0$
  of the vertex of the reentrant corner.
  For $\varphi = \varphi^*$, we have the asymptotic  relation
\[
\theta^* (\rho) = \theta^* + \frac{2 (a_1 b_2 - a_2 b_1)\,
\tan \frac{\varphi^*}{\beta}}{(a_1^2 + b_1^2)^{3/2}}\,
\,\rho \,+\,\mathcal{O} (\rho^{2}),
\qquad \rho \to 0,
\]
 where $\theta^* (\rho)$ denotes $\theta_{\varphi^*} (\rho)$.
 For other values of $\varphi$, we have the asymptotic  relation
\[
\theta_\varphi (\rho) = \theta^* +
 \frac{\mu b_1 \sin (\varphi - \varphi^*)}
{(a_1^2 + b_1^2)^{\frac{\beta + 1}{2}} \cos^\beta \frac{\varphi}{\beta}}
\rho^{\beta -1} +
 \mathcal{O} (\rho^{2 \beta - 2}), \quad \rho \to 0, \quad
 \varphi^* \not= \varphi \in \Big(- \frac{\pi\beta}{2},\, \frac{\pi\beta}{2}\Big).
\]

\begin{figure}
\centering
\resizebox*{12cm}{!}{\includegraphics{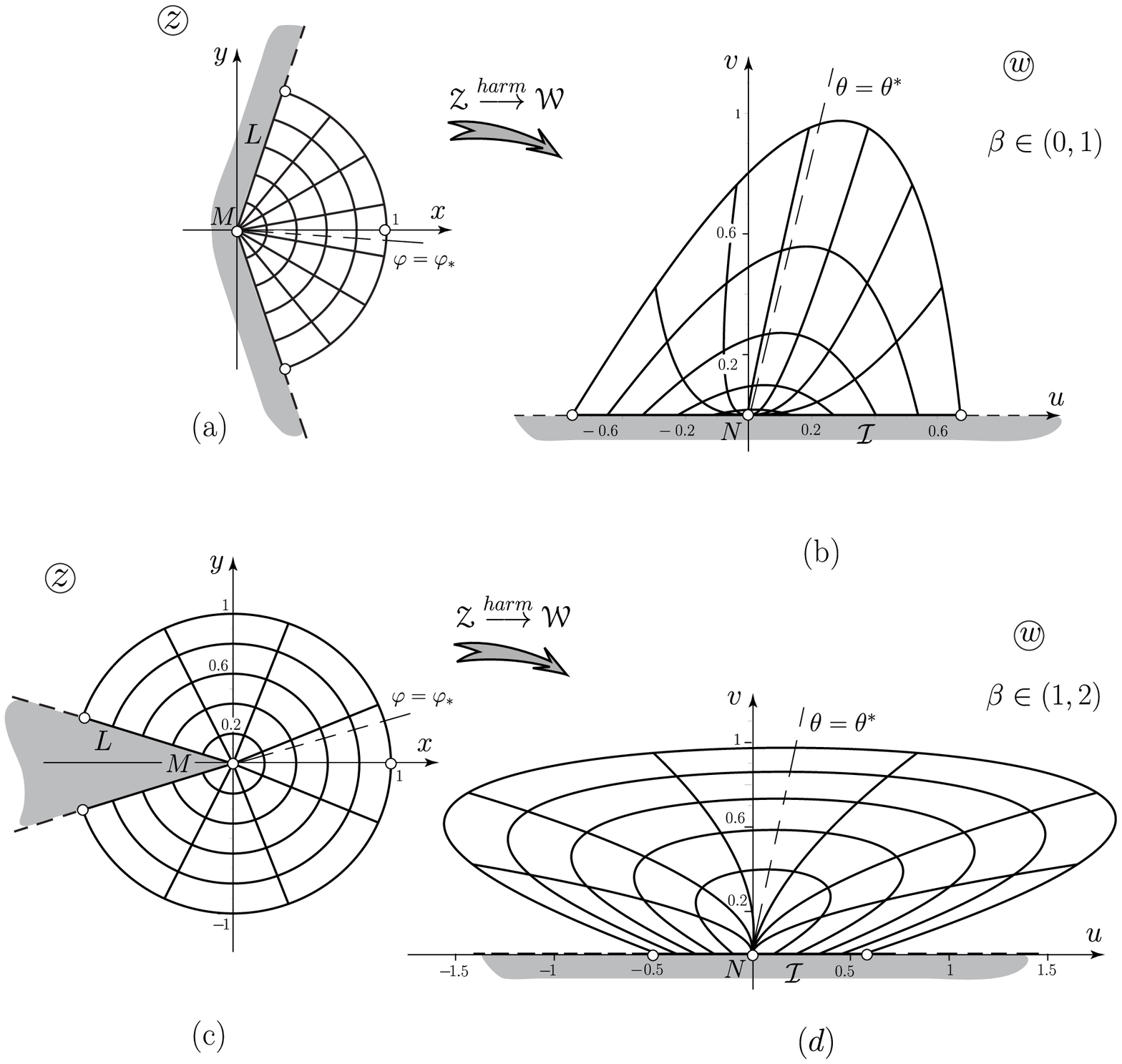}}
\caption{} \label{fig5}
\end{figure}

Combining the results of the present section, we arrive at the
following assertion.

\begin{theorem}\label{theo-2}
    If  $\beta \in (0, 1),$ then the asymptotic behavior of the function $v (u)$
    describing the curve $l_\varphi = \mathcal{F} (\lambda_\varphi)$ in the  Cartesian
    coordinates is expressed by relations \eqref{BV-18b}-\eqref{BV-18c}.
 If $\beta \in (1, 2),$ then the specified curve is described by
 means of asymptotics  \eqref{BV-20}-\eqref{BV-21}
   in the  Cartesian
    coordinates
$(\mathcal{U}, \mathcal{V})$.
\end{theorem}

It follows that the dependence $\theta (\varphi)$ between the angle $\varphi$, 
under which the ray $\lambda_\varphi$ leaves the vertex $z = 0$ of the corner, 
and the angle $\theta$, under which its image 
$l_\varphi: = \mathcal {F}^{-l} (\lambda_\varphi)$ leaves the preimage of 
the vertex $w = 0$, is described for $\beta \in (0, 1)$ by the formula
\begin{equation}
\beta \in (0, 1): \qquad \theta (\varphi) \,=\, \left\{\begin{array}{cl} 0, 
&\varphi \in [- \pi \beta / 2,\, \varphi^*), \\
\; \theta^*, &\varphi = \varphi^*, \\
\pi, &\varphi \in (\varphi^*, \pi \beta / 2].
\end{array}\right., 
\label{BV-22}
\end{equation}
and when $\beta \in (1, 2)$ --- by the following formula:
\begin{equation}
\beta \in (1, 2): \qquad \theta (\varphi) \,=\, \left\{\begin{array}{cl} 0, 
&\varphi = - \pi \beta / 2, \\ 
\; \theta^*, &\varphi \in (- \pi \beta / 2,\, \pi \beta / 2), \\
\pi , &\phi = \pi \beta / 2.
\end{array}\right. 
\label{BV-23}
\end{equation}
In Fig. 5a is a polar grid
\[
T = \Big\{r = \frac{n}{5}, 
n = \overline{1, 5};\;\varphi = \frac{\pi \beta}{2} \Big(\frac{2n}{7} - 1\Big),
n = \overline{0, 7}\Big\},
\]  
and in Fig. 5b and 5c show its images $\mathcal{F} (T)$ on the plane $w$, respectively 
for $\beta \in (0, 1)$ and $\beta \in (1, 2)$.

Theorems \ref{theo-1}, \ref{theo-2} admit the following obvious adherence.

\begin{theorem}
Let the boundary of the domain $\mathcal{Z}$ contain the contour of the angle
$L = L_- \cup M \cup L_+$ of the solution $\pi \beta \in (0, 2 \pi)$ with vertex
$M$, $\mathcal{W}$ contains rectilinear segment 
$\mathcal{I} = \mathcal{I}_- \cup N \cup \mathcal{I}_+$, and the map
$\mathcal{F}: \mathcal{Z} {\buildrel {\, harm \,} \over \longrightarrow} \mathcal {W}$, 
univalent and preserving orientation near $L$, corresponds to the following conditions:
$\mathcal{F} (L) = \mathcal{I}$, $\mathcal{F} (M) = N$, and for the uniform motion
of a point $z \in L_\mp$ the speed of its image on 
$\mathcal{I}_\mp$ is equal to the constant $\sigma_\mp
$.

Further, let $\mathcal{Z} \supset \lambda_\varphi$ be a smooth Jordan curve,
leaving the vertex $M$ at an angle $\varphi$, $\theta$ is the angle at which its image
leaves the point $N$, and $\theta (\varphi)$ is the correspondence between these angles; 
let it go $\mathcal{W} \supset \Lambda_\varphi$ is a smooth Jordan curve issuing 
from the point $N$ at an angle $\theta$, $\varphi$ is the angle by which its image
$\mathcal{F}^{-1} (\Lambda_\theta)$ leaves the vertex $M$, and $\varphi (\theta)$
--- the correspondence between these angles.

Then the dependence $\theta (\varphi)$ for $\beta \in (0, 1)$
is given by the formula (\ref{BV-22}), and for $\beta \in (1, 2)$ 
by the formula (\ref{BV-23}); the dependence $\varphi (\theta)$
for $\beta \in (0, 1)$ is given by the formula (\ref{BV-14}), and for 
$\beta \in (1, 2)$, by the formula (\ref{BV-18a}).
\end{theorem}

Thus, the dependencies $\theta (\varphi)$ and $\varphi (\theta) $ are discontinuous
for a harmonic mapping of a domain with an angle, in contrast to a similar conformal
imaging
$\mathcal{K}: \mathcal{S} {\,\buildrel {conf\,} \over\longrightarrow\,} \mathcal{H}$, 
$\mathcal{K} (0) = 0$, for which these dependences are linear,
$\theta (\varphi) = \varphi / \beta + \pi / 2$, 
$\varphi (\theta) = \beta \theta - \pi \beta / 2$. 

\vspace*{3mm}

{\bf Acknowledgments} 
 This work was supported by Russian Foundation for Basic Research,
 projects
 16--01--00781
  and
   16--07--01195,
    and by the Ministry of Education
 and Science of the Russian Federation on the program to improve
 the competitiveness of Peoples' Friendship University (RUDN University)
  among the world's leading
 research and education centers in the 2016--2020.


\begin{thebibliography}{10}

\bibitem{Hamilton1975}
Hamilton, R. Harmonic maps of manifolds with boundary,
 Berlin--Heidelberg--N-Y.:
  Springer--Verlag,
   1975.




\bibitem{Jost1985}
Jost, J. Harmonic mappings and minimal immersions,
 Berlin--N-Y.:
  Springer--Verlag,
   1985.


\bibitem{Duren2004}
 Duren, P. Harmonic mappings in the plane,
  Cambrige: Cambrige University Press, 2004.



\bibitem{Bshouty2005}
 Bshouty, D., Hengartner, W. Univalent harmonic mappings in the plane,
 Handbook of complex analysis: geometric function
theory, Vol. 2, Amsterdam: Elsevier, 2005, pp. 479--506.




\bibitem{Rado1926}
Rad\'o, T. Aufgabe 41, Jahresbericht Deutsch. Math.--Verein.,
 1926, Vol. 35, p. 49.

\bibitem{Kneser1926}
Kneser, H. L\"osung der Aufgabe 41, Jahresbericht Deutsch.
Math.--Verein.,
 1926, Vol. 35, pp.  123--124.

\bibitem{Choquet1945}
    Choquet, G. Sur un type de transformation analitiques
g\'en\'eralisant la repr\'esentation conforme et d\'efinie au
moyen de fonctions harmoniques, Bull. Sci. Math., 1945,
 Vol. 69, No 2, pp.  156--165.


\bibitem{BezrodnykhVlasov-2012}
 Bezrodnykh, S. I., Vlasov, V. I. On a problem in the constructive theory
  of harmonic mappings, J. Math. Sci. (N.Y.), 2014,
 Vol. 201, No 6, pp. 705--732.

\bibitem{Caratheodory1913}
 Caratheodory, C. \"Uber die gegenseitige Bezeihung der R\"ander bei der
Konformer Abbildung des Inneren einer Jordanschen Kurve auf einer
Kreis, Math. Ann., 1913,
 Vol. 73, pp. 305--320.



\bibitem{Renelt1988}
 Renelt, H. Elliptic systems and quasiconformal mappings,
 N.-Y.: John Wiley \& Sons, 1988.

\bibitem{Liao1991}
 Liao, G. On harmonic maps,
 Mathematical aspects of numerical grid
generation, Philadelphia: SIAM, 1991, pp. 123--130.

\bibitem{Alessandrini2009}
 Alessandrini, G., Nesi, V. Invertible harmonic mappings, beyong Kneser,
 Ann. Scuola Norm, Sup. Pisa, Cl. Sci., 2009,
 Vol. 8, No 3, pp. 451--468.



\bibitem{Bers1951}
 Bers, L. Isolated singularities of minimal surfaces,
 Ann. of Math., 1951,
 Vol. 53,  364--386.


\bibitem{Nitsche1958}
 Nitsche, J. C. C. On an estimate for the curvature of minimal surfaces $z =
z(x, y)^*$,
 Journal of Mathematics and Mechanics, 1958,
 Vol. 7, No 5, pp. 767--769.


\bibitem{Clunie1984}
 Clunie, J,
 Sheil--Small, T. Harmonic univalent functions,
 Ann. Acad. Sci. Fenn. Math., 1984,
 Vol. 9,   pp. 3--25.

\bibitem{Hall1985}
 Hall, R. A class of isoperimetric inequalities,
 J. Analyse Math., 1985,
 Vol. 45,   pp. 169--180.



\bibitem{SheilSmall1990}
 Sheil--Small, T. Constants for planar harmonic mappings,
 J. London Math. Soc., 1990,
 Vol. 42,   pp. 237--248.


\bibitem{KlajPonnusamyVuorinen2013}
Klaj, D., Ponnusamy, S., Vuorinen, M. Radius of close-to-convexity
and fully starlikeness of harmonic mappings,
 Complex Var. Elliptic Equ., 2013,
 Vol. 59, No 4, pp. 539--552.

\bibitem{BshoutyLundbergWeitsman2015}
 Bshouty, D., Lundberg, E., Weitsman, A. A solution to Sheil--Small's
  harmonic mapping
problem for plygons,
 Proc. Amer. Math. Soc., 2015,
 Vol. 143, No 12, pp. 5219--5225.

\bibitem{Bagapsh2017}
 Bagapsh, A. O. On the radius of starlikeness for harmonic mappings,
 Zap. Nauchn. Sem. S.-Peterburg.
  Otdel. Mat. Inst. Steklov. (POMI), 2014,
 Vol. 456, pp. 16--24.



\bibitem{Grotzsch1928}
 Gr\"otzsch, H. \"Uber die Verzerrung bei schlichten nichtkonformen
Abbildungen und \"uber eine damit zusammenh\"angende Erweiterung
des Picardschen Satzes,
 Berichte Leipzig, 1928,
 Vol. 80,   pp. 503--507.

\bibitem{Morrey1938}
 Morrey, Ch. B. Jr. On the solutions of quasi-linear
  elliptic partial
differential equations,
 Trans. Amer. Math. Soc., 1938,
 Vol.  43, No 1 , pp. 126--166.

\bibitem{Teichmuller1940}
 Teichm\"uller, O. Extremal quasikonforme Abbildungen und
quadratische Differentiale,
 Abh. Preuss. Akad. Wiss. Math.-Nat.
 Kl, 1940,
 Vol. 1939, No 22, pp. 3--197.

\bibitem{Lavrentiev1948}
 Lavrentiev, M. A. The basic theorem of the theory of quasi--conformal
mappings of planar domains,
 Izv. Akad. Nauk SSSR, Ser. Mat., 1948,
 Vol. 12, No 6, pp. 513--554.

\bibitem{Hersch1952}
 Hersch, J., Pfluger, A. G\'en\'eralisation du lemme de Schwarz et du principe
de la mesure harmonique pour les fonctions pseudo-analytiques,
 C. R. Math. Acad. Sci. Paris, 1952,
 Vol. 234,   pp.  43--45.

\bibitem{Ahlfors1966}
 Ahlfors, L. V. Lectures on quasiconformal mappings,
 N.-Y.--London:
 D. Van Nostrand Co., 1966.

\bibitem{Martio1968}
 Martio, O. On harmonic quasiconformal mappings,
 Ann. Acad. Sci. Fenn. Math., 1969,
 No 425,   pp. 1--10.



\bibitem{Lehto1973}
 Lehto, O., Virtanen, K. I. Quasiconformal mappings in the plane,
  Berlin--Heidelberg--N-Y.:
  Springer--Verlag, 1973.

\bibitem{Vuorinen1997}
 Anderson, G. D.,  Vamanamurthy, M. K., Vuorinen M. Conformal invariants,
  inequalities, and quasiconformal maps,
  Chichester: Wiley, 1997.


\bibitem{Bers1953}
Bers, L. Univalent solutions of linear elliptic systems,
 Comm. Pure Appl. Math., 1953,
 Vol. 6,  pp. 513--526.

\bibitem{BojarskiIwaniec1974}
Bojarski, B. V., Iwaniec, T. Quasiconformal mappings and
non-linear elliptic equations in two variables,
 Bull. Pol. Acad. Sci. Math., 1974,
 Vol. 22,  pp.  473--484.

\bibitem{HengartnerSchober1986-5}
 Hengartner, W., Schober, G. Harmonic mappings with given dilatation,
 J. London Math. Soc., 1986,
 Vol. 33,   pp.  473--483.


\bibitem{DurenKhavinson1997}
 Duren, P., Khavinson, D. Boundary correspondence and dilatation of harmonic
mappings,
 Complex Var. Theory Appl., 1997,
 Vol. 33,   pp.   105--111.



\bibitem{Winslow1967}
 Winslow, A. Numerical solution of the quazi-linear
  Poisson equations in a nonuniform triangle mesh,
 J. Comput. Phys., 1967,
 Vol. 2,   pp.  149--172.

\bibitem{Godunov1972}
 Godunov, S. K., Prokopov, G. P.
 The utilization of movable grids in gas dynamic calculations,
 Zh. Vychisl. Mat. Mat. Fiz., 1972,
 Vol. 12, No 2,  pp. 429-440.


\bibitem{Brackbill1982}
 Brackbill, J. U., Saltzman, J. S. Adaptive zoning for singular problems in two dimensions,
 J. Comput. Phys., 1982,
 Vol. 46, No 3,   pp. 342--368.

\bibitem{Roache1985}
 Roache, P. J., Steinberg, S. A new approach to grid genaration using a variational
formulation,
 Proc. of AIAA 7th Computational Fluid Dynamics Conference, Cincinnati, 1985, pp. 360--370.

\bibitem{WarsiTuarn1986}
 Warsi, Z. U., Tuarn, W. N. Surface mesh generation using elliptic equations,
 Numerical grid generation in computational fluid
dynamics, Whiting, NJ: Pineridge Press, 1986, pp.  95--100.

\bibitem{Sengupta1988}
 Sengupta, S. et al.(Eds.) Numerical grid generation in computational fluid
mechanics, Whiting, NJ: Pineridge Press, 1989.


\bibitem{KnSt1993}
 Knupp, P., Steinberg, S. Fundamentals of grid generation,
 Boca Raton, FL: CRC Press, 1993.


\bibitem{Ivanenko1997}
 Ivanenko, S. A. Adaptive harmonic meshes,
 M.: Computing Center of the RAS, 1997.

\bibitem{Smith1998}
 Smith, P. W., Sritharan, S. S. Theory of harmonic grid generation,
 Complex Var. Theory Appl., 1988,
 Vol. 10, No 4,  pp.  359--369.

\bibitem{ThSW1999}
 Thompson, J. F., Soni, B. K., Weatherill, N. P. (Eds). Handbook of grid generation,
 Boca Raton, FL: CRC Press, 1999.

\bibitem{Liseikin1999}
 Liseikin, V. D. Grid generation methods. scientific computation,
 N.-Y.: Springer-Verlag,
 1999.


\bibitem{Azarenok2014}
 Azarenok, B. N., Charakhch'yan, A. A. One problem on generating
  2D regular grids based on mappings,
 Math. Models Comput. Simul., 2015,
 Vol. 7, No 4,  pp. 303--314.

\end{thebibliography}
\end{document}